# On $L_1$-distance between first exit times from two regions[*]

Nikolai Dokuchaev

Department of Mathematics and Computer Sciences,

The University of West Indies, Mona, Kingston 7, Jamaica W.I.

October 30, 2018

## Abstract

First exit times from regions and their dependence on variations of boundaries are discussed for diffusion processes. The paper presents an estimate of $L_1$-distance between exit times from two regions via expectations of exit times.

*Abbreviated title:* On $L_1$-distance between first exit times

It is known that first exit times from a region for smooth solutions of ordinary equations does not depend continuously on variations of the initial data or on the boundary of the region. But first exit times for non-smooth trajectories of diffusion processes have some path-wise regularity with respect to these variations (some related results can be found in author's papers (1987),(1992)). This paper studies path-wise dependence on fluctuations of the boundary for first exit times of diffusion processes. We present an effective estimate of $L_1$-distance between exit times from two regions for a diffusion process via expectations of exit times.

Let $(\Omega, \mathcal{F}, \mathbf{P})$ be a standard probability space. Consider a $n$-dimensional diffusion process $y(t)$ such that

$$dy(t) = f(y(t))dt + \beta(y(t))dw(t), \quad t > 0, \tag{1}$$

$$y(0) = a. \tag{2}$$

Here $w(t)$ is a standard $d$-dimensional Wiener process, $f$ and $\beta$ are non-random functions with respective values in $\mathbf{R}^n$ and $\mathbf{R}^{n \times d}$. The random vector $a$ with values in $\mathbf{R}^n$ does not depend on $w(\cdot)$.

All vectors and matrices are real, and $\bar{D}$ denotes the closure of a region $D$.

---

Further, for $x \in \mathbf{R}^n$, we denote by $y^x(t)$ the solution of (1) with the initial condition $y(0) = x$.

For a set $\Gamma \in \mathbf{R}^n$, we denote $\tau(\Gamma) \triangleq \inf\{t: y(t) \in \Gamma\}$. This is the first times of achieving $D$ of the processes $y(t)$. Similarly, we denote $\tau^x(\Gamma) \triangleq \inf\{t: y^x(t) \in \Gamma\}$.

Let $D_1$, $D_2$ be two bounded regions in $\mathbf{R}^n$. Let $\Gamma_i = \partial D_i$ be the boundary of $D_i$.

We assume that the boundaries of $D_i$ are $C^1$-smooth. Further, we assume that all the components of the functions $f$, $\beta$ are continuously differentiable, and $b(x)b(x)^\top \geq cI_n$, where $I_n$ is the unit matrix in $\mathbf{R}^{n \times n}$, and $c > 0$ is a constant.

**Theorem 1** *Let $a \in \bar{D}_1 \cap \bar{D}_2$ with probability 1. Then*

$$\mathbf{E}|\tau(\Gamma_1) - \tau(\Gamma_2)| \leq \max\left(\sup_{x \in D_1 \cap \Gamma_2} \mathbf{E}\tau^x(\Gamma_1), \sup_{x \in D_2 \cap \Gamma_1} \mathbf{E}\tau^x(\Gamma_2)\right). \tag{3}$$

Note the theorem is oriented on the case when $D_1 \setminus D_2 \neq \emptyset$ and $D_2 \setminus D_1 \neq \emptyset$. If $D_1 \subset D_2$ or $D_2 \subset D_1$ then the estimation (3) is obvious.

**Example.** Let $n = d = 1$, $y^a(t) = a + w(t)$, $D_0 = (0,1)$, $D_\varepsilon = (\varepsilon, 1+\varepsilon)$, where $\varepsilon \in [0,1)$. Let $\Gamma_\varepsilon \triangleq \{\varepsilon, 1+\varepsilon\}$. Let $\tau_\varepsilon^a = \inf\{t: y^a(t) \in \Gamma_\varepsilon\}$. We have that $\mathbf{E}\tau_\varepsilon^x = (x-\varepsilon)(1+\varepsilon-x)$ for any $x \in D_\varepsilon$ (it can be found similarly (4) and (9) below). By Theorem 1, it follows that

$$\mathbf{E}|\tau_0^a - \tau_\varepsilon^a| \leq \max\left(\mathbf{E}\tau_\varepsilon^0, \mathbf{E}\tau_0^\varepsilon\right) = \varepsilon(1-\varepsilon).$$

**Remark** We assumed that the boundaries and coefficients are smooth, the diffusion is non-degenerate, and the regions are bounded. In fact, these conditions can be lifted provided that the right hand part of (3) is finite. For example, a similar theorem can be obtained for first exit times of the process $(y(t), t)$ from cylindrical regions $D_i \times (0, T)$, $i = 1, 2$, $T > 0$, for the case of time dependent coefficients $f$ and $\beta$.

*Proof of Theorem 1.* Let $v_i = v_i(x) : D_i \to \mathbf{R}$, $i = 1, 2$, be the solutions in $D_i$ of the Dirichlet problems

$$\mathcal{L}v_i = -1, \quad v_i|_{\Gamma_i} = 0. \tag{4}$$

Here the differential operator

$$\mathcal{L} = \sum_{k=1}^n f_k \frac{\partial}{\partial y_k} + \frac{1}{2} \sum_{k,l=1}^n b_{k,l} \frac{\partial^2}{\partial y_k \partial y_l}, \tag{5}$$

where $f_k$, $y_l$, and $b_{k,l}$ are the components of the vectors $f$, $y$ and the matrix $b = \beta\beta^\top$. As is known, the problem (4) has the unique solution that is twice continuously differentiable up to the boundary.

Let $e_1$ and $e_2$ be the indicator functions of the random events $\{\tau(\Gamma_1) > \tau(\Gamma_2)\}$ and $\{\tau(\Gamma_2) > \tau(\Gamma_1)\}$ respectively.

Let $\mathcal{F}_t$ be the filtration generated by $w(t)$ and $a$.

Let $\widehat{\tau} \triangleq \tau(\Gamma_1) \wedge \tau(\Gamma_2)$. The random variables $e_i$ are measurable with respect to the $\sigma$-algebras $\mathcal{F}_{\widehat{\tau}}$, $\mathcal{F}_{\tau(\Gamma_i)}$, $i = 1, 2$, associated with the Markov times (with respect to the filtration $\mathcal{F}_t$) $\widehat{\tau}$ and $\tau(\Gamma_i)$ (see, e.g., Gihman and Skorohod (1975), Chap. 4, §2). Using Itô's formula, we obtain the equality

$$\begin{aligned}
\mathbf{E}\left\{e_1 v_1[y(\tau(\Gamma_2))]\right\} &= -\mathbf{E}\left\{e_1\{v_1[y(\tau(\Gamma_1))] - v_1[y(\tau(\Gamma_2))]\}\right\} \\
&= -\mathbf{E}\left\{e_1 \int_{\widehat{\tau}}^{\tau(\Gamma_i)} \mathcal{L}_1 v_1[y(t)] dt\right\} \\
&= \mathbf{E}\left\{e_1[\tau(\Gamma_1) - \widehat{\tau}]\right\} \\
&= \mathbf{E}\left\{e_1[\tau(\Gamma_1) - \tau(\Gamma_2)]\right\}.
\end{aligned} \quad (6)$$

If we replaced the indices $1, 2$ in (6) by 2,1, we get similarly that

$$\mathbf{E}\left\{e_2 v_2[y(\tau(\Gamma_1))]\right\} = \mathbf{E}\left\{e_2[\tau(\Gamma_2) - \tau(\Gamma_1)]\right\}. \quad (7)$$

Clearly,

$$\mathbf{E}|\tau(\Gamma_1) - \tau(\Gamma_2)| = \mathbf{E}\left\{e_1[\tau(\Gamma_1) - \tau(\Gamma_2)]\right\} + \mathbf{E}\left\{e_2[\tau(\Gamma_2) - \tau(\Gamma_1)]\right\}. \quad (8)$$

We have that

$$v_i(x) = \mathbf{E}\, \tau^x(\Gamma_i). \quad (9)$$

Then it follows from (6)-(8) that

$$\begin{aligned}
\mathbf{E}|\tau(\Gamma_1) - \tau(\Gamma_2)| &\leq \mathbf{E}\left\{e_1[\tau(\Gamma_1) - \tau(\Gamma_2)]\right\} + \mathbf{E}\left\{e_2[\tau(\Gamma_2) - \tau(\Gamma_1)]\right\} \\
&= \mathbf{E}\left\{e_1\{v_1[y(\tau(\Gamma_2))]\}\right\} + \mathbf{E}\left\{e_2\{v_2[y(\tau(\Gamma_1))]\}\right\} \\
&\leq \max\left(\sup_{x \in D_1 \cap \Gamma_2} v_1(x), \sup_{x \in D_2 \cap \Gamma_1} v_2(x)\right).
\end{aligned}$$

Now the assertion of the Theorem 1 follows. $\square$

Nikolai Dokuchaev,

Department of Mathematics and Computer Science, The University of West Indies,

Mona, Kingston 7, Jamaica W.I.

email ndokuch@uwimona.edu.jm